\documentclass[11pt,a4paper]{article}
\usepackage[utf8]{inputenc}
\usepackage{amsmath, amssymb, amsthm}
\usepackage{graphicx,hyperref,geometry}
\usepackage{a4wide}
\usepackage{xcolor,verbatim}

\newtheorem{theorem}[subsection]{Theorem}%[section]
\newtheorem{corollary}[subsection]{Corollary}
\newtheorem{lemma}[subsection]{Lemma}

\newtheorem{conjecture}[subsection]{Conjecture}

\newtheorem{question}[subsection]{Question}
\newtheorem{problem}[subsection]{Problem}

\newcommand{\cH}{\mathcal{H}}
\newcommand{\cA}{\mathcal{A}}
\newcommand{\cB}{\mathcal{B}}

\title{The connection between the chromatic numbers of a hypergraph and its $1$-intersection graph}

\author{Zoltán L. Blázsik\thanks{\protect\includegraphics[height=1cm]{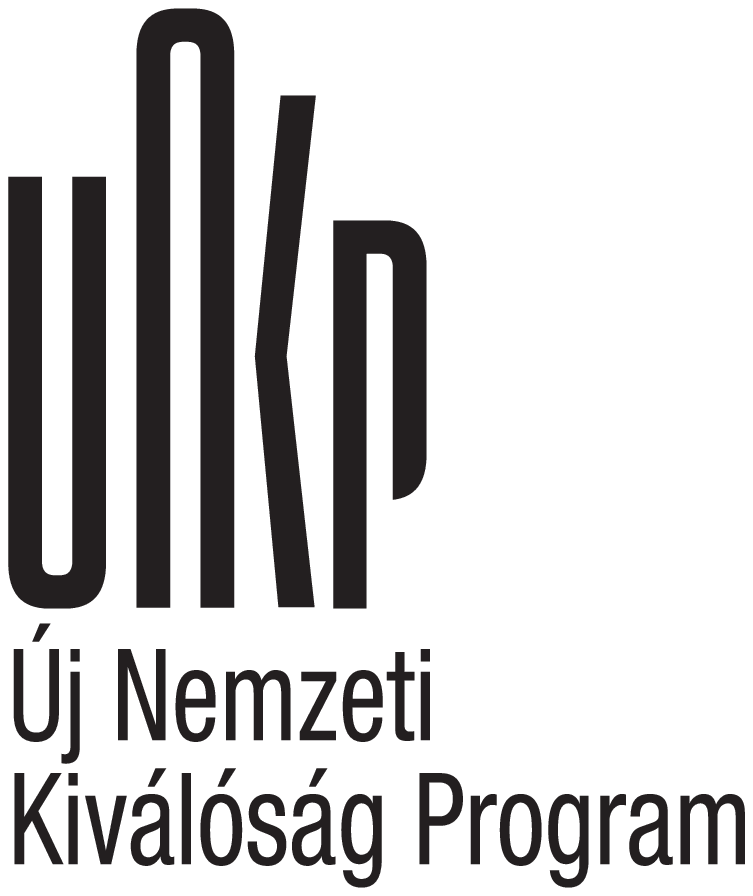}\includegraphics[height=0.8cm]{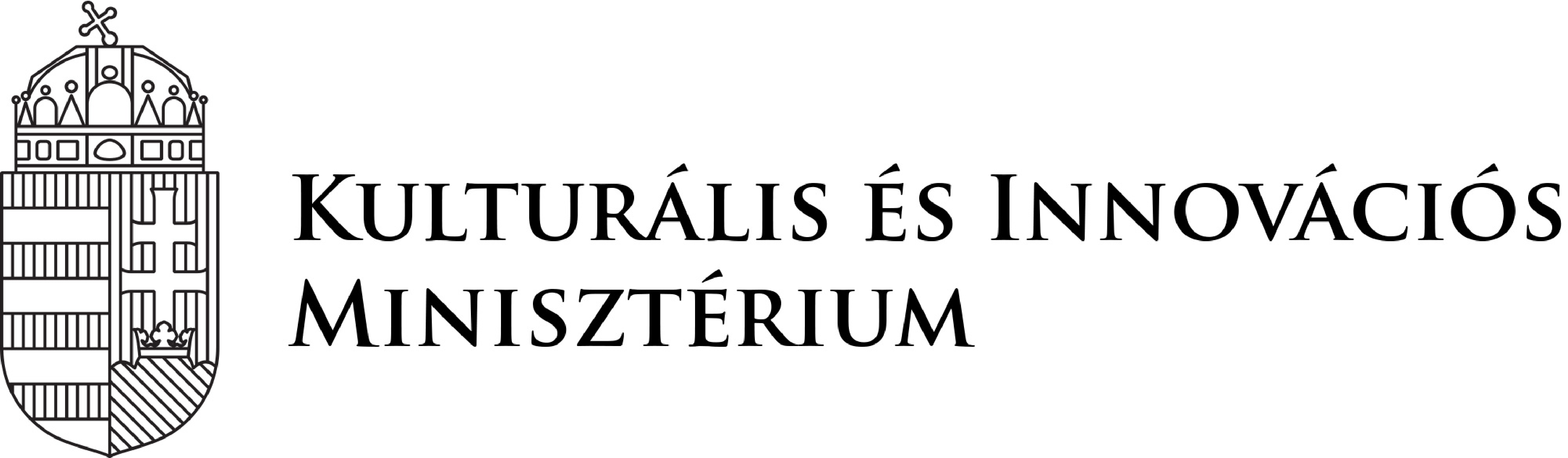} The author was supported by the ÚNKP-23-4-SZTE-628 New National Excellence Program of the Ministry for Culture and Innovation from the source of the National Research, Development and Innovation Fund, and also by the OTKA grant no. SNN 132625.} \\
\small HUN-REN Alfr\'ed R\'enyi Institute of Mathematics\\
\small SZTE Bolyai Institute \\
\small \url{blazsik@server.math.u-szeged.hu} \\
and\\
Nathan W. Lemons \\ %\thanks{~{\color{red}Please fill it with your grants.}}
\small Theoretical Division\\
\small Los Alamos National Laboratory\\
\small \url{nlemons@lanl.gov} \\
}
\date{\today}

\begin{document}

\maketitle

\begin{abstract}
A well known problem from an excellent book of Lovász states that any hypergraph with the property that no pair of hyperedges intersect in exactly one vertex can be properly 2-colored. Motivated by this as well as recent works of Keszegh and of Gyárfás et al we study the $1$-intersection graph of a hypergraph. The $1$-intersection graph encodes those pairs of hyperedges in a hypergraph that intersect in exactly one vertex. We prove for $k\in\{2,4\}$ that all hypergraphs whose $1$-intersection graph is $k$-partite can be properly $k$-colored.
\end{abstract}

\section{Introduction}
Coloring combinatorial objects to achieve certain properties is a major theme of combinatorial theory.  Graphs, and their generalization hypergraphs, are defined by the pair $V,E$ where $E$ is a collection of subsets of $V$.  If the edges (elements of $E$) all have size $r$ we call the object an $r$-graph or an $r$-uniform hypergraph; 2-graphs are better known as simply graphs. A coloring problem of particular interest is the study of which hypergraphs and graphs can be vertex-colored with $c$ colors
so that no edges are monochromatic. In the sequel, we assume that any hyperedge has size at least $2$ which guarantees that there always exists a vertex coloring without monochromatic hyperedges. Such vertex colorings are often called proper.

This problem has led to many important and well-known results.  Graphs that are 2-colorable, that is bipartite, are well characterized and easy to determine.  On the other hand, it is well known that determining if a graph is 3-colorable (also called 3-chromatic) is NP-hard\cite{Karp}.  For hypergraphs, it is already NP-hard to determine if a given hypergraph can be 2-colored.  Unsurprisingly, there are many results on the (relative) hardness of coloring fixed subclasses of graphs and hypergraphs.  For instance, the famous four-color theorem states that any planar graph can be properly colored with 4 colors~\cite{AH4}. The interested reader is pointed to the surveys \cite{CroftJensen,MolloyReed,Lewis} for more information on coloring graphs.

A well known problem (Problem 13.33.) from the excellent book of Lov\'asz states that any hypergraph with the property that no pair of edges intersect in exactly one vertex can be 2-colored \cite{Lovász}.  For two edges in a hypergraph, we say edge $e$  has a 1-intersection if there exists another edge $f$ with $|e\cap f|=1$.  The inductive proof of this elegant problem from Lov\'asz' book actually implies more: that if one can properly 2-color those edges of a hypergraph that have 1-intersections, then the whole hypergraph can be properly two colored.  Thus it is natural to focus on those edges with 1-intersections.  In recent work, Keszegh did this for hypergraphs whose edges with 1-intersections satisfied a certain local, geometric condition \cite{Keszegh}.

We consider a different, global condition on the 1-intersections of a hypergraph.  For a hypergraph $\mathcal{H}$ we define the 1-intersection graph $\cH^{[1]}$ of $\mathcal{H}$ to be a graph on the edge set of $\mathcal{H}$.  That is the vertices of this graph will be $E(\mathcal{H})$.  Two hyperedges $e,f\in E(\mathcal{H})$ are adjacent in $\mathcal{H}^{[1]}$ if and only if $|e\cap f|=1$.
Gy\'arf\'as et al \cite{Gyarfas} asked the following question.
\begin{question}\label{prob-uniform}
Let $r\geq 3$ and let $\mathcal{H}$ be an $r$-uniform hypergraph whose 1-intersection graph $\cH^{[1]}$ is $k$-colorable for some $k\geq 2$. Then is $\mathcal{H}$ also $k$-colorable? In other words is $\chi(\cH)\le \chi(\cH^{[1]})=k$ true?
\end{question}
The main result of their paper was giving a positive answer for the $r=3$ case of this question.

\begin{theorem}[Gy\'arf\'as et al]
   Let $\mathcal{H}$ be a 3-uniform hypergraph and suppose $\chi(\cH^{[1]})=k\ge 2$. Then $\chi(\mathcal{H})\le \chi(\mathcal{H}^{[1]})=k$.
\end{theorem}
In this paper, we focus on a natural generalization of Question \ref{prob-uniform} that concerns non-uniform hypergraphs:

\begin{question}\label{prob-nonuniform}
    Let $\mathcal{H}$ be a hypergraph such that $\chi(\mathcal{H}^{[1]})=k\ge 2$. Is $\mathcal{H}$ $k$-colorable?
\end{question}
We remark that the greedy coloring argument of the motivational problem automatically shows that if the 1-intersection graph of $\cH$ is $k$-colorable, then $\mathcal{H}$ is $(k+1)$-colorable. Our main result is to give a positive answer to Question \ref{prob-nonuniform} for the case when $k\in\{2,4\}$.

\begin{theorem}\label{t:main}
For any hypergraph $\cH$ if $\chi(\cH^{[1]})=2$ then $\chi(\cH)=2$.
\end{theorem}

\begin{theorem}\label{t:main4}
For any hypergraph $\cH$ if $\chi(\cH^{[1]})=4$ then $\chi(\cH)\le4$.
\end{theorem}

Let us remark here that the inequality is important because it might happen that the chromatic number of the original hypergraph $\cH$ is even smaller. We note that because it is easy (algorithmically) to tell if a graph is bipartite, our results imply that it is easy algorithmically to show certain hypergraphs are 2-colorable. Namely one has only to show that the 1-intersection graph is bipartite. On the other hand, the converse is not necessarily true. Consider a hypergraph such that all the hyperedges have a common universal vertex. Then the chromatic number of the 1-intersection graph is equal to the number of edges since they form a clique in the $1$-intersection graph. Since every hyperedge has size at least 2 hence the hypergraph is 2-colorable because we get a proper $2$-coloring by giving color red to the universal vertex and blue to every other vertex.  On the other hand there are examples where the chromatic number of the hypergraph is equal to the chromatic number of its $1$-intersection graph.  Namely, for all $\ell\ge 1$ the complete graph $K_{2\ell+1}$ satisfies this.  For the even case, we can take the complete graph $K_{2\ell}$ with the union of one arbitrary $3$-edge.

Finally we note that Question \ref{prob-nonuniform} has a negative answer in the following cases.  Let $\ell\geq 1$ and consider the graph $K_{2\ell}.$  The edges can be written as the disjoint union of $2\ell-1$ perfect matchings.  Thus the 1-intersection graph is $(2\ell-1)$-partite, however one needs $2\ell$ colors to properly color $K_{2\ell}$.  By Brooks' theorem~\cite{Brooks} these are the only 2-graphs which can give a negative answer to Question \ref{prob-nonuniform}.  It would be interesting to know if there are any hypergraphs which give a negative answer for Question \ref{prob-nonuniform} with $k$ even. Is there any hypergraphs which give a negative answer for Question \ref{prob-nonuniform} such that every hyperedge has size at least 3?

The paper is organised as follows. In Section 2, we first prove Theorem \ref{t:main} using a finite recoloring process and then we prove Theorem \ref{t:main4} using our previous result for $k=2$. In Section 3, we finish with some open problems.

\section{Main results}

To prove our main result for $k=2$, we use induction on the number of edges in a hypergraph with bipartite $1$-intersection graph.  This is achieved by employing a vertex recoloring process which maintains a queue of currently monochromatic hyperedges and in each step it corrects the next monochromatic hyperedge by recoloring a well-chosen vertex.

\begin{proof}[Proof of Theorem \ref{t:main}]
The proof is by induction on $|E(\cH)|$; it is easy to see that the base case is satisfied.

Suppose now that the statement in Theorem \ref{t:main} holds for any hypergraph with at most $m$ edges. Let $\cH$ be a hypergraph with $m+1$ edges such that the 1-intersection graph is bipartite.  Let $h_0$ be an edge of $\cH$.  By the induction hypothesis, $\cH-{h_0}$ can be properly 2-colored.  Now consider this same coloring on $\cH$.  Without loss of generality we can assume that $h_0$ is monochromatic in this coloring.

We define a process to recolor some of the vertices of $\cH$.  At the end of the recoloring process, no edges will be monochromatic.
The process works through a queue of currently monochromatic edges, denoted by the ordered tuple $Q_t$ at step $t$.  The process starts with $Q_0=(h_0)$ and at each further step $t$, let $h_t$ denote the first edge in the tuple $Q_t$.  If there is some vertex $v_t\in h_t$ which has not yet been recolored, then (exactly one such vertex) $v_t$ is recolored to a new color. Note that there is only one choice of color.  This new coloring may introduce further monochromatic edges, denote these by $C_t$.  However, $h_t$ is of course no longer monochromatic.  The queue $Q_{t+1}$ is obtained from $Q_t$ by removing $h_t$ (and any other edges of the queue that are no longer monochromatic) and by adding the newly monochromatic edges $C_t$.  If any edges of $C_t$ have already been in the queue at any previous time, the process stops with failure.  Likewise, if every vertex of $h_t$ has already been recolored the process stops with failure.
The proof will be complete when we show that the process does not fail but successfully recolors the vertices of $\cH$.

    Consider the bipartite $1$-intersection graph of $\cH$ and denote the two parts by $\cA$ and $\cB$. Suppose that $h_0\in \cA$ is monochromatic blue; denote the other color by red.  Clearly every edge of $\cH$ belongs to exactly one class $\cA$ or $\cB$.

    \begin{lemma}\label{l:recoloring}
        During the recoloring process, any hyperedge which is (at any time during the process) monochromatic blue must belong to $\cA$.  Likewise, any monochromatic red hyperedge must belong to $\cB$.
    \end{lemma}
\begin{proof}
    We will show this is true by induction on $t$, the steps of the process.  Clearly at $t=0$, there is only one monochromatic edge: it is $h_0$ and is both blue and in $\cA$.  Now suppose for some $t$ all the monochromatic blue edges of $\cup_{i=0}^t Q_i$ are in $\cA$ and all the monochromatic red edges of $\cup_{i=0}^t Q_i$ are in $\cB$.  We will show that the monochromatic blue and monochromatic red edges of $Q_{t+1}$ must also belong to $\cA$ and $\cB$ respectively.  Suppose $h_t$ is monochromatic red.  Then by hypothesis, $h_t$ is a member of $\cB$.  In addition, during the $t^{\textrm{th}}$ step, a vertex of $h_t$ is recolored from red to blue.  Thus any monochromatic edges of $Q_{t+1}\backslash\big(\cup_{i=0}^t Q_i\big)$ must be monochromatic blue.  But these edges  intersect $h_t$ in exactly one vertex; that is they must all belong to $\cA$.  An identical argument holds if $h_t$ is monochromatic blue: in this case all the newly monochromatic edges will be red and belong to $\cB$.
\end{proof}

\begin{corollary}\label{c:1}
    No edge can be added twice to the queue.
\end{corollary}
\begin{proof}
   Suppose an edge $h$ is added twice.  As the process stops before recoloring any vertex twice, $h$ must have started the process monochromatic, that is $h$ can only be $h_0$.  But the second time $h_0$ is added to the queue it must be monochromatic red.  But then $h_0$ must belong to both $\cA$ and $\cB$, a contradiction.
\end{proof}

\begin{corollary}\label{c:2}
    At each time step $t$ of the process, the edge $h_t$ must have at least one vertex which has not yet been recolored.
\end{corollary}
\begin{proof}
    Suppose not.  Let $t$ be a time step of the process such that all the vertices of $h_t$ have already been recolored once.  Then all the vertices of $h_t$ have been recolored exactly once as no vertex is recolored more than once.  But then $h_t$ must be the edge $h_0$ as all other edges of the hypergraph start the process properly colored.  But this is a contradiction as $h_0$ has then been added to the queue a second time.
\end{proof}

By definition the process is finite.  Corollaries \ref{c:1} and \ref{c:2} imply that the process does not fail.  Thus it must terminate by reaching a proper coloring of $\cH$.
\end{proof}

The proof of our second main result for the $k=4$ case uses the result of Theorem \ref{t:main} as a black box, but unfortunately the similar argument will not work for $k=2\ell$ if $\ell\ge 3$.

\begin{proof}[Proof of Theorem \ref{t:main4}]
Consider a proper $4$-coloring $c^{[1]}$ of the $1$-intersection graph $\cH^{[1]}$. Let us denote the color classes with $C_1$, $C_2$, $C_3$ and $C_4$. Consider two subhypergraphs of $\cH$, denoted by $\cH_{12}$ and $\cH_{34}$, both of which are spanning (i.e. the vertex set is $V$) and the set of hyperedges of $\cH$ are partitioned into two parts: $\cH_{12}$ contains those hyperedges of $\cH$ that belong to $C_1 \cup C_2$, and $\cH_{34}$ contains those hyperedges of $\cH$ that belong to $C_3 \cup C_4$, respectively.

Observe that the $1$-intersection graphs $\cH_{12}^{[1]}$ and $\cH_{34}^{[1]}$ of the two subhypergraphs are bipartite (or equivalently $2$-colorable) by definition. By Theorem \ref{t:main}, there are two, possibly different, proper $2$-colorings $c_{12}$ and $c_{34}$ of $V$ such that $c_{12}$ with colors $0$ and $1$ is a proper $2$-coloring of $\cH_{12}$, and $c_{34}$ with colors $0'$ and $2$ is a proper $2$-coloring of $\cH_{34}$. Now we define a new $4$-coloring $c$ of $V$ by by taking the ,,sum'' of $c_{12}$ and $c_{34}$ in the following way.

\[
(\forall v\in V) \qquad c(v)=\left \{\begin{tabular}{ll}
   0 & $\quad$ if $c_{12}(v)=0$ and $c_{34}(v)=0'$, \\
   1 & $\quad$ if $c_{12}(v)=1$ and $c_{34}(v)=0'$, \\
   2 & $\quad$ if $c_{12}(v)=0$ and $c_{34}(v)=2$, \\
   3 & $\quad$ if $c_{12}(v)=1$ and $c_{34}(v)=2$.
\end{tabular}\right .
\]

We claim that $c$ is a proper $4$-coloring of $\cH$. Suppose on the contrary that there exists a hyperedge $h=\{v_1,v_2,\dots,v_r\}$ which is monochromatic with respect to $c$. No matter which partition class $C_i$ contains $h$, if $c(v_1)=c(v_2)=\dots=c(v_r)$ then $c_{12}(v_1)=c_{12}(v_2)=\dots=c_{12}(v_r)$ and $c_{34}(v_1)=c_{34}(v_2)=\dots=c_{34}(v_r)$ must hold, too. Equivalently, since $h$ is monochromatic with respect to $c$ it must have been monochromatic with respect to $c_{12}$ or $c_{34}$ depending on the color class of $h$ with respect to $c^{[1]}$, which is a contradiction. Hence $c$ is a proper $4$-coloring of $\cH$, therefore $\chi(\cH)\le 4$.
\end{proof}

\section{Concluding remarks}

It is very intriguing that in the general case the implication's correctness depends heavily on the parity of the chromatic number of the $1$-intersection graph. We believe that the statement should hold for any $k$ even values as well.

\begin{conjecture}\label{c:main_even}
For any hypergraph $\cH$ if $\chi(\cH^{[1]})=2\ell$ then $\chi(\cH)\le2\ell$.
\end{conjecture}

Let us conclude our paper by stating some open problems.

\begin{problem}
    Is there any counterexample if $\chi(\cH^{[1]})=2\ell$ for some $\ell\ge2$?
\end{problem}

\begin{problem}
    Is there any counterexample if the minimum size of the hyperedges of $\cH$ is at least 3?
\end{problem}

\section*{Acknowledgement}
This research was started during the 11th Emléktábla Workshop, 2022. The authors are thankful for the organizers for inviting them, and especially to Viola Mészáros, Cory Palmer for useful discussions on the topic which was introduced to us by Balázs Keszegh.

\end{document}